\documentclass[11pt]{amsart}
\usepackage{amssymb}
\usepackage{palatino}
\usepackage{tikz-cd}
\input amssym.def

\usepackage{amsmath, amsfonts}
\usepackage{amssymb}
\usepackage{amscd}
\usepackage[mathscr]{eucal}
\usepackage{palatino}
\usepackage{listings}

\newfont{\cyrr}{wncyr10}
\usepackage{tikz-cd}
\newcommand{\thmref}[1]{Theorem~\ref{#1}}
\newcommand{\corref}[1]{Corollary~\ref{#1}}

\newcommand{\lemref}[1]{Lemma~\ref{#1}}
\newcommand{\K}{{\mathbf K}}
\renewcommand{\H}{{\mathbf H}}
\newcommand{\F}{{\mathbf F}}
\newcommand{\Cl}{{\rm Cl}}
\renewcommand{\b}{{\mathfrak{b}}}
\renewcommand{\a}{{\mathfrak{a}}}

\renewcommand{\P}{{\mathfrak{p}}}

\renewcommand{\O}{{{\mathcal{O}}}}

\newcommand{\N}{{\mathbb N}}
\newcommand{\Z}{{\mathbb Z}}
\newcommand{\Q}{{\mathbb Q}}

\newcommand{\G}{{\mathbf{G}}}

\newtheorem{thm}{Theorem}

\newtheorem{lem}[thm]{Lemma}
\newtheorem{cor}[thm]{Corollary}

\newtheorem{rmk}{Remark}[section] 
\newtheorem{defn}{Definition}

\title{Euclidean ideal classes in Galois number fields of odd prime degree}
\author{V. Kumar Murty and J. Sivaraman}

\address[ V. Kumar Murty ]
{Department of Mathematics,
University of Toronto,
40 St. George Street,
Toronto, ON, Canada, M5S 2E4.
}

\email{murty@math.toronto.edu} 
\address[J. Sivaraman]   
{Chennai Mathematical Institute,
H1, SIPCOT IT Park, Siruseri,
Kelambakkam, India, 603103. 
}

\email{jyothsnaas@cmi.ac.in}


\begin{document}

\maketitle

\begin{abstract}
Weinberger \cite{PW} in 1972, proved that the ring of integers
of a number field with unit rank at least $1$ is a principal ideal domain if and
only if it is a Euclidean domain, provided the
generalised Riemann hypothesis holds. Lenstra \cite{HL},
extended the notion of Euclidean domains in
order to capture Dedekind domains with finite
cyclic class group and proved an analogous
theorem in this setup. More precisely, he showed
that the class group of the ring of integers of a number field with
unit rank at least $1$ is cyclic if and only
if it has a Euclidean ideal class, provided the
generalised Riemann hypothesis holds. 
The aim of this paper is to show the
following.
Suppose that $\K_1$
and $\K_2$ are two Galois number fields of odd prime degree with cyclic 
class groups and Hilbert class fields
that are abelian over $\Q$. 
If $\K_1\K_2$ is ramified
over $\K_i$, 
then at least one $\K_i$ ($i \in \{1,2\}$)
must have a Euclidean ideal class.
\end{abstract}
\section{Introduction}
A well known result of Weinberger \cite{PW} from 1972
states that the ring of integers of a number field
$\K$ with unit rank at least $1$,
has class number $1$ if and only if it
is a Euclidean domain, provided we assume
that the generalised Riemann hypothesis
(Riemann hypothesis for all Dedekind zeta functions) 
holds.
In 1977, this result was extended by
 Lenstra \cite{HL-EA}. He showed that 
if the ring of S-integers of a number
field with $|S| \ge 2$ is a principal
ideal domain, it is Euclidean,
again provided we assume the
generalised Reimann hypothesis.

The first unconditional result
on this problem came in a 
paper of the first author, Gupta and Ram Murty \cite{GMM}
in 1985.
Using the techniques introduced in
the famous work of Gupta and Ram Murty on
Artin's conjecture \cite{GM}, 
the authors of \cite{GMM} were able to
show unconditionally that rings of S-integers
in certain Galois number fields $\K$
have trivial S-class group if and only
if they are Euclidean
when $|S| \ge \max\{5, 2[\K:\Q]-3\}$.

In 1995,
Clark and Ram Murty \cite{CM}
were able to extend 
this result to rings of integers of totally
real Galois number fields of
degree atleast $4$ with an additional
property. At this juncture,
the case of lower degree had
still not been addressed.

In 2004, Harper in his thesis,
showed that the ring $\Z[\sqrt{14}]$
is a Euclidean domain \cite{MH}.
This was immediately followed by
a paper of Harper and Ram Murty who introduced 
some important new ideas 
to show unconditionally
that the ring of integers of a Galois number
field with unit rank atleast 4 is Euclidean if 
and only if it has class number $1$ \cite{HM}. 
Further in the case of abelian
number fields they were able to extend
the same result to the case of unit rank $3$.
Narkiewicz \cite{WN}
in 2007, used some of the ideas from
the works of Harper \cite{MH}, Harper and Ram Murty~\cite{HM}  
and Heath-Brown \cite{RH}
to address
the case of lower unit rank in greater detail. 
He showed that for
Galois cubic fields with class number~$1$,
there is atmost $1$ whose ring of integers is not a Euclidean
domain and in case of real quadratic fields
with class number $1$ there are atmost
$2$ whose rings of integers are not Euclidean domains.

On the other hand, in 1979, Lenstra \cite{HL} extended the notion
of Euclidean domains to cover a larger family of
Dedekind domains with finite cyclic class groups,
thereby extending the notion beyond class number
$1$. In order to do so he defined Euclidean ideal classes. 
We state the definition of Euclidean ideal classes
below while assuming that the ambient ring
is the ring of integers of a number field.
\begin{defn}
	
	Let $\O_{\K}$ be the ring of integers of a number field $\K$, $E$ be the set of all non-zero integral 
	ideals of $\O_{\K}$ and $W$ be a well-ordered set.  Given a map 
	$$
	\psi : E \to W,
	$$
	we say that a non-zero ideal $\a$ of $\O_{\K}$ is Euclidean for $\psi$ if 
	for each non-zero 
	ideal $\b$ of $\O_{\K}$ and any $x~\in~{\b}^{-1}{\a} - \a$, 
	there exists $y \in \a$ such that
	$$
	\psi\left( {\a}^{-1}{\b(x-y)} \right) < \psi(\b).
	$$
	Further the class $[\a]$ is called a Euclidean
	ideal class of $\O_{\K}$.
\end{defn}

Lenstra's arguments show that if $\O_{\K}$ has a Euclidean ideal class, it must have
a finite cyclic class group which is in fact generated
by the Euclidean ideal class. He goes on to show
that for number fields with unit rank atleast $1$,
if the class group is cyclic it must have a Euclidean
ideal class, provided the generalised Riemann
hypothesis holds.

The first unconditional result in the
Euclidean ideal class setup was shown by
Graves and Ram Murty \cite{GM}. Their precise
result is stated below.
\begin{thm}{\rm(Graves and Ram Murty \cite{GM})}
	Let $\K$ be a number field with ring of integers 
	$\mathcal{O}_{\K}$ and cyclic class group $Cl_{\K}$. If its Hilbert class
	field $\H(\K)$, has an abelian Galois group over $\Q$ and if the unit 
	rank of $\mathcal{O}_{\K}^{\times}$ is atleast $4$, 
	then  the class group $\Cl_{\K} = \langle[C]\rangle$ if and only if $[C]$ is a Euclidean 
	ideal class.
\end{thm} 

The result was extended to a family of number
fields with unit rank $3$ by Deshouillers, Gun
and Sivaraman \cite{DGS}. We state the precise
result below.
	\begin{thm}{\rm (Deshouillers, Gun and Sivaraman~\cite{DGS})}
	Suppose that $\K$ is a number field with 
	unit rank greater than or equal to $3$
	and its Hilbert class field $\H(\K)$ is abelian over $\Q$.
	Also suppose that the conductor of $\K$ (that is, the least positive integer $n$ for which $\K \subset \Q(\zeta_n)$) is $f$ and $\Q(\zeta_f)$ 
	over $\K$ is cyclic. Then the class group 
	$\Cl_{\K}$ is cyclic if and only if it has 
	a Euclidean ideal class. 
\end{thm} 
	
Investigations have also been carried out
in number fields of lower unit rank.
The theorem of 
Narkiewicz \cite{WN} was extended to this setup
under certain conditions by Gun and Sivaraman~\cite{GS}. We state below the result of Gun and
Sivaraman \cite{GS} in the case of Galois cubic fields.

Let $\K_1$ and $\K_2$ be number fields 
with Hilbert class fields $\H(\K_1)$ and  $\H(\K_2)$
(abelian over $\Q$) respectively. 
Also let $f_1$ and $f_2$ be the conductors of $\H(\K_1)$ and 
$\H(\K_2)$ respectively. 
Set $f$ to be the least common multiple of $16, f_1$
and $f_2$. Further,
$\F := \Q(\zeta_{f})$, where $\zeta_f$ is a primitive $f$th
root of unity.  
In this set up, the result is as follows.

\begin{thm}{\rm(Gun and Sivaraman~\cite{GS})}
	Let $\K_1, \K_2$ be distinct Galois cubic fields 
	with prime class numbers and $\H(\K_1), \H(\K_2), {\F}$, $f$ 
	be as above.  Also let $G$ be the Galois group of 	
	$\F$ over $\K_{1}\K_{2}$, 
	$G_{\ell}$ be the Galois group of $\F$ over $\Q(\zeta_{\ell})$,
	where either $\ell$ is an odd prime dividing $f$ or $\ell=4$. If
	$$
	G \not\subset \bigcup_{\ell} G_{\ell} \bigcup \text{ Gal} (\F/\H(\K_1))
	\bigcup \text{ Gal} (\F/\H(\K_2)),
	$$
	then at least one of $\K_1, \K_2$ has 
	a Euclidean ideal class.  
\end{thm}

In this paper we extend the above result 
by removing the condition $G \not\subset \bigcup_{\ell} G_{\ell}$
where $\ell$ is either an odd prime dividing $f$ or $\ell = 4$.
The precise result is stated below.
\begin{thm}\label{main-cubic}
	Let $\K_1$ and $\K_2$ be distinct Galois number fields of odd prime
	degree
	with cyclic class groups,
	such that the Hilbert class fields are abelian over $\Q$. 
	If $\K_1\K_2$ is ramified over $\K_i$,
	then at least one $\K_i$ for $i \in \{ 1,2 \}$ has 
	a Euclidean ideal class. 
\end{thm}

The rest of the paper is organised as follows.
In the preliminaries, we state some results
which will be used in due course of the main
proof. This will be followed by a section on
choosing a residue class modulo
$f$. In the next section, we give the proof of Theorem \ref{main-cubic}
using this choice of residue class. In the penultimate section, we give
some examples of fields for which our theorem holds and the final section
consists of the data availability statement.
 \section{preliminaries}
 
 We state below some definitions and theorems 
 required to complete the proof of
 our theorem.
 We first state the definition
 of a genus field of an abelian
 number field $\K$.
 \begin{defn}
 For an abelian number field $\K$,
 the maximal extension $\G(\K)$ of $\K$, abelian over $\Q$,
 which is unramified at all
 the finite places of $\K$ is called
 the genus field of $\K$. Further
 $[\G(\K):\K]$ is called the genus number
 of $\K$, denoted $g_{\K}$.
 \end{defn}
 The first theorem we state here
 is due to Leopoldt stated as a corollary
 in Section 5 of~\cite{HLeo}.
 
 \begin{thm}{\rm(Leopoldt \cite{HLeo}, Pg 52 \cite{MI})}~\label{GT}
 For an abelian number field $\K$ the genus
 number $g_{\K}$ is given by
 $$
 g_{\K} = \frac{\prod_q e(q)}{[\K:\Q]}
 $$
 where the product runs over rational primes $q$ and $e(q)$ is the ramification index of $q$ with
 respect to the field $\K$.
 \end{thm}
 
 We now introduce the definition
 of the Hilbert class field of a number
 field $\K$.
 
 \begin{defn}\label{HCF}
 For a number field $\K$, the Hilbert
 class field $\H(\K)$ is the maximal extension
 of $\K$, abelian over $\K$, which is unramified at all the places
 of $\K$. The degree $[\H(\K) : \K]$
is denoted by $h_{\K}$.
 \end{defn}
 
 \begin{rmk}
 We note that if the Hilbert class field $\H(\K)$
 is abelian over $\Q$, then it is an abelian
 number field containing $\K$ and unramified
 at all the finite places of $\K$. Therefore
 $\H(\K) \subset \G(\K)$ and therefore
 $h_{\K} \mid g_{\K}$. In particular
 for a field $\K$ of odd prime degree $p$ whose $\H(\K)$
 is abelian over $\Q$, the ramification indices $e(q) \mid p$ for all rational primes $q$.
 Therefore $h_{\K}$ is a non-negative power of $p$ (by \thmref{GT}).
 This fact will be used crucially in the proof of \lemref{generator}.
 \end{rmk}
 
 The next theorem is 
 Theorem 14 of Gun and Sivaraman \cite{GS},
 on Euclidean ideal classes in
 quadratic and cubic fields.
 \medskip
 
 \begin{thm}{\rm(Gun and Sivaraman \cite{GS})}~\label{HarpCritvar}
 	Suppose that $\K$ is a number field with unit rank 
 	greater than or equal to one and its class group
 	$\Cl_{\K} = \langle[\a]\rangle$. 
 	If there exists an unbounded increasing
 	sequence $\{x_{n}\}_{n \in \N}$ such that  
 	\begin{equation*}
 	|\{\mathfrak{p} ~:~\P \text{ is prime},  ~[\mathfrak{p}]= [\a],  
 	\mathfrak{N}(\P) \le x_n, 
 	\text{ every residue class of }(\mathcal{O}_{\K}/\mathfrak{p})^{\times} 
 	\text{ contains a unit }\}| 
 	~\gg~ 
 	\frac{x_n}{ \log^{2} x_n}~,
 	\end{equation*}
 	then $[\a]$ is a Euclidean ideal class of $\O_{\K}$.
 \end{thm}
 The next is Lemma 3 of Heath-Brown's 
 paper on Artin's conjecture \cite{RH}.
 \begin{lem}{\rm (Heath-Brown~\cite{RH})}~\label{Heathbrown}
 	Suppose that $u$ and $v$ are natural numbers with the
 	following properties
 	$$
 	(u,v)=1, 
 	\phantom{m} 
 	v\equiv 0 \bmod16
 	\phantom{m}\text{and}\phantom{m}
 	\left(\frac{u-1}{2}, ~v\right) ~=~1.
 	$$ 
 	Then there exist $a, b \in (\frac{1}{4}, \frac{1}{2})$ 
 	with $a <b$ such that for any $\epsilon > 0$, the 
 	set
 	\begin{eqnarray*}
 		P(X) 
 		& : = & 
 		\{p\equiv u \bmod v ~:~ p \in (X^{1-\epsilon}, X) \text{ such that } 
 		\frac{p-1}{2} \text{ is either prime or }\\ && 
 		~\phantom{m}\text{ is a product of primes }q_{1}q_{2} 
 		\text{ with } X^{a}\le q_{1} \le X^{b}\}
		 \end{eqnarray*}
		has cardinality $\gg X/\log^2 X$.
 \end{lem}
 
 The final lemma is Lemma 2 of 
 Narkiewicz's paper \cite{WN}, on  counting units in residue
 classes of number fields.\\
 \begin{lem}[{\rm Narkiewicz~\cite{WN}}]~\label{primroot}
 	Let $\K$ be an arbitrary algebraic number field and
 	let $a_{1},a_{2}$ and $a_{3}$ be multiplicatively independent 
 	elements of $\K^{\times}$, $T$ be a 
 	set of prime ideals of degree $1$ in $\K$
 	and $\mathfrak{N}_{\K/\Q}(\P)$ denotes the absolute
 	norm of a prime ideal $\P$  of $\mathcal{O}_{\K}$.
 	Suppose that $T$ has the following properties;
 	\begin{enumerate}
 		\item 
 		there exists a constant $c > 0$ and an unbounded 
 		increasing sequence $\{ x_n \}_{n \in \N}$ such that
 		$$
 		|T(x_n) := \{ \P \in T ~:~ \mathfrak{N}_{\K/\Q}(\P) \le x_n \}| 
 		~>~ cx_n/\log^{2}x_n   ~~\text{ for all  } n.
 		$$
 		\item 
 		there exist $\alpha, \beta \in (1/4, 1/2)$ 
 		with $\alpha < \beta$ such that if $\mathfrak{p} \in T$ 
 		and $p := \mathfrak{N}_{\K/\Q}(\mathfrak{p})$, then either 
 		$p-1= 2q$ or $p-1 = 2q_{1}q_{2}$ where
 		$q$, $q_1$ and $q_2$ are primes
 		such that $p^{\alpha} < q_1 < p^{\beta}$.
 		\item 
 		the numbers $a_1, a_2$ and $a_3$ are 
 		quadratic non-residues with respect to 
 		every prime in $T$.
 	\end{enumerate}
 	Then for any $0 < \epsilon < c$, there exists a 
 	subsequence $\{ y_m\}_{m \in \N}$ of 
 	$\{ x_n \}_{n \in \N}$ such that 
 	one of the $a_{i}$s is a primitive root 
 	for at least $(c-\epsilon)y_m /\log^2y_m$ 
 	elements of $T(y_m)$.
 \end{lem}
In the next section, we will construct
an appropriate residue class before
we proceed to the proof of the main
theorem. 
\section{Choosing a residue class}
We briefly state some notation that 
will be used in this section.
Let $\K_1$ and $\K_2$ be distinct fields of degree $p_1$ and $p_2$ (odd primes), with
cyclic class groups and Hilbert class fields $\H(\K_1)$ and  $\H(\K_2)$
(assumed abelian over $\Q$) respectively. 
Also let $f_1$ and $f_2$ be the conductors of $\H(\K_1)$ and 
$\H(\K_2)$ respectively. 
Set $f$ to be the least common multiple of $16, f_1$
and $f_2$. Further let
$\F := \Q(\zeta_{f})$, where $\zeta_f$ is a primitive $f$th
root of unity. Let $G$ be the Galois group
of $\F$ over $\K_1\K_2$. 
\begin{lem}\label{generator}
	If $\K_1\K_2$ is ramified over $\K_i$, then
	there exists a  residue class $d \bmod f \in G$ such that
	the image of $ d \bmod f $ in $Gal(\H(\K_i)/\K_i)$
	is a generator of the group and $d \bmod f \notin 
	Gal(\F/\Q(\zeta_\ell))$ for any odd prime $\ell \mid f$ or $\ell =4$.
	
\end{lem}

\begin{proof}
	
If $h_{\K_i} \neq 1$ and $\K_1\K_2$ is ramified
over $\K_i$, it
follows that
$\H(\K_i)\cap \K_1\K_2 = \K_i$.
If $h_{\K_i}=1$ we still have
$\H(\K_i)\cap \K_1\K_2 = \K_i$.
Therefore we have the following diagram
when $\H(\K_i) \neq \K_i$ for $i \in \{1,2\}$.
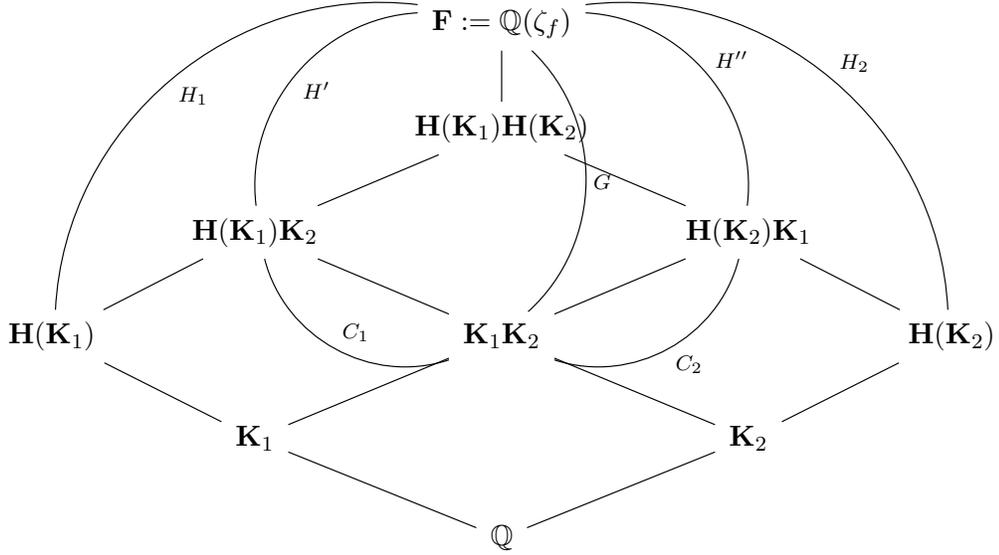
$$
\begin{tikzcd}
&
&
&
\F:=\Q(\zeta_{f})
\arrow[bend left=47, dash]{ddd}{G}
\arrow[bend right=47, dash]{dddll}{H_1}
\arrow[bend right=47, dash]{ddl}{H'}
\arrow[bend left=47, dash]{ddr}{H''}
\arrow[bend left=47, dash]{dddrr}{H_2}
\arrow[dash]{d} \\
&
&
&
\H(\K_1)\H(\K_2) \arrow[dash]{dr}{}
\arrow[dash]{dl}{}\\
&
&
\H(\K_1) \K_2
 \arrow[dash]{dr} 
\arrow[dash]{dl}{}
\arrow[bend right=47, dash]{dr}{C_1}
&
&  
\H(\K_2) \K_1 \arrow[dash]{dl}
\arrow[bend left=47, dash]{dl}{C_2}
\arrow[dash]{dr}\\
&
\H(\K_1) \arrow[dash]{dr}
&
&
\K_1\K_2 \arrow[dash]{dl}
\arrow[dash]{dr}
&
&
\H(\K_2) \arrow[dash]{dl}\\
&
&
\K_1 \arrow[dash]{dr}{}
& 
&
\K_2 \arrow[dash]{dl}{}\\
&
&
&
\Q
\end{tikzcd}
$$

Consider a $c \bmod f \in G$
which restricts to a generator of $C_1$ (Galois group of $\H(\K_1)\K_2/\K_1\K_2)$ and 
$C_2$ (Galois group of $\H(\K_2)\K_1/\K_1\K_2$). Note that these groups
are non-trivial and isomorphic to the class groups of $\K_1$ and $\K_2$, respectively and cyclic by our assumptions.
This also implies that $1  \not\equiv c \bmod f$.\\

\textbf{Claim:} We may assume that the order of
$c \bmod f$ is odd.\\
If not, we may replace $c \bmod f$ with
$c' \equiv c^{2^m} \bmod f$ where the $2$-adic valuation of the order of $c \bmod f$
 is $m$.
Such a $c' \bmod f$ is still in $G$ and since the
class numbers are powers of  the odd primes $p_1$ and $p_2$ (by the remarks after Definition \ref{HCF}), it will still restrict to a generator
of $C_1$ and $C_2$.\\

Let $H_i$ be the Galois group of $\Q(\zeta_f)/H(\K_i)$, $H'$
the Galois group of $ \Q(\zeta_f)/H(\K_1)\K_2$
and $H''$ the Galois group of $ \Q(\zeta_f)/H(\K_2)\K_1$.
Further, since the $\K_i$ are totally real, $-1 \bmod f \in H_1 
\cap H_2 \cap G = H' \cap H''$.
This implies that $\overline{c}H' = -\overline{c}H'$ and $\overline{c}H'' = -\overline{c}H''$.
So we consider the element $-c \bmod f \in G$
such that
\begin{enumerate}
	\item $c \bmod f $ has odd order;
	\item $c \bmod f $ restricts to a generator
	of $C_1$ and $C_2$ and therefore so does
	$-c \bmod f$.
\end{enumerate}

We observe that $ (-c \bmod f)^{ord(c \bmod f)}  \equiv -1 \bmod f$. This implies that 
$$-1 \bmod f \in \langle-c \bmod f\rangle.$$
This implies that $-c \bmod f \notin Gal(\F/\Q(\zeta_\ell))$ for any 
odd prime $\ell \mid f$ or $\ell=4$. In other words
$\left(\frac{-c-1}{2}, f\right)~=~1$.

If only one of the class groups is trivial
then we can find a $c \bmod f$
which corresponds to a generator
of the non-trivial class group with the
above two properties. 

The resulting congruence
class $-c \bmod f$ is in $G$ and therefore 
corresponds to a generator
for both class groups.

If both the class groups are trivial then
$c \equiv -1\bmod f$ will 
satisfy the required conditions.
\end{proof}
 
\smallskip
 \section{Proof of \thmref{main-cubic}}

Throughout this subsection, let
$\K_1$ and $\K_2$ be Galois fields of odd prime degree 
with Hilbert class fields $\H(\K_1)$ and 
$\H(\K_2)$, abelian over $\Q$ respectively. 
Also let $f_1$ and $f_2$ be their conductors.
Set $f$ to be the least common multiple of $16, f_1, f_2$.  
Choose a positive integer $u_1$ such that  it satisfies $u_1 \equiv d \bmod f$.
It follows from \lemref{Heathbrown} and \lemref{generator} that there exist $a, b \in (\frac{1}{4}, \frac{1}{2})$ 
with $a <b$ such that for any $\epsilon >0$
\begin{eqnarray*}
|J_{\epsilon}(X)
&:=&
\{p\equiv u_1 \bmod f ~:~ p \text{ rational prime, } 
p \in (X^{1-\epsilon}, X) \text{ such that }\\
&& 
\phantom{mm}\frac{p-1}{2} \text{ is either a rational}
 \text{  prime or a product of rational} \\
&& 
\phantom{mm} \text{primes }
q_{1}q_{2} 
\text{ with } X^{a}<q_{1}< X^{b}\} | 
\gg 
\frac{X}{\log^2 X} .
\end{eqnarray*}

Next set $\K=\K_{1}\K_{2}$. Note that any prime 
$p \in \Z$ which is congruent to $u_1 \bmod f$
splits in $\K_{1}\K_{2}$. Since $\left(\frac{u_1-1}{2},f\right)=1$, 
we note that $u_1 \equiv 3 \bmod 4$.
Choose $\epsilon$ such that $a < \frac{b}{1-\epsilon} < \frac{1}{2}$.
Consider the set
\begin{eqnarray*}
M_{\epsilon}
&:=&
\{\mathfrak{p}~:~ \P \text{ is a prime ideal in } \O_{\K}, ~ \mathfrak{N}_{\K/\Q} (\mathfrak{p})
= p \text{ rational prime }, p \equiv u_1 \bmod f,\\
&&\phantom{mmmm}  
~\frac{p-1}{2} \text{ is either a rational prime or a product of } \phantom{m}\\  
&&\phantom{mmmm} \text{  rational 
primes } q_{1}q_{2}  \text{  with } 
p^{a} < q_{1} < p^{\frac{b}{1-\epsilon}}\};
\end{eqnarray*}
and also the set 
$M_{\epsilon}(X)~:=~\{\P \in M_{\epsilon}: \mathfrak{N}_{\K/\Q}(\P)\le X\}$
for any real number $X > 0$. 
It follows that $M_{\epsilon}(X) \gg \frac{X}{\log^2 X}$.
Suppose that $\epsilon_{1},\epsilon_{2}$ and
$\epsilon_{3}$ are fundamental units from $\K_i$
with atleast one from each $\K_i$ for $i \in \{1,2\}$.   \\
\\
Following the proof of Lemma 16 in \cite{GS},
we write $M_{\epsilon} = \cup_{n=1}^{8} M_{n}$.
where each $M_n$ corresponds to a tuple
$(c_1,c_2,c_3)$ with entries in ${\pm 1}$
such that $ \left(\frac{\epsilon_i}{\P}\right) = -c_i$ for
all $\P \in M_n$. Here $\left(\frac{\cdot}{\P}\right)$ 
is used to denote the second power residue symbol.\\
\\
Since $M_{\epsilon}(X) \gg \frac{X}{\log^2 X}$,
we claim that there exists an increasing 
unbounded sequence of positive real numbers $\{x_m\}_{m \in \N}$ and a $c >0$ such that 
$$M_{n_0}(x_m) > c\frac{x_m}{\log^2 x_m} ~~\text{ for some } ~1 \le n_0 \le 8.$$
Suppose not,
then for any $1 \le n \le 8$, we have
\begin{eqnarray*}
\limsup_{X \to \infty} \frac{M_n(X)}{X/\log^2 X} = 0 & \text{ and } &
\liminf_{X \to \infty} \frac{M_n(X)}{X/\log^2 X} = 0.
\end{eqnarray*}
This implies $M_{\epsilon}(X) = o\left(\frac{X}{\log^2 X}\right)$, which
is a contradiction. 

Further since $u_1 \equiv 3 \bmod 4$, it follows that $\left(\frac{-1}{\P}\right)=-1$. 
For the tuple $(c_1,c_2,c_3)$ corresponding to $M_{n_0}$, $c_i\epsilon_i$ 
is a quadratic non-residue modulo any prime $\P \in M_{n_0}$. Now,
on applying \lemref{primroot} with $T= M_{n_0}$ there exists
a subsequence $\{y_l\}_{l \in \N}$ of $\{x_m\}_{m \in \N}$
such that one of the
elements $\pm \epsilon_1, \pm \epsilon_{2}, \pm \epsilon_{3}$ is a primitive root modulo $\P$
for atleast $(c-\epsilon)y_l/\log^2 y_l$ elements
of $M_{n_0}(y_l)$ for any $0 < \epsilon < c$. 
We denote this subset of elements of $M_{n_0}$ by $Q$ and the primitive root by $\eta$.
\medskip

 For $Q$ 
and $\eta$ as above, it follows that $\eta \in \K_{s}$ for some 
$s \in \{1,2\}$. Further we have a sequence 
$\{ y_l \}_{l \in \N}$ such that 
\begin{eqnarray*}
Q(y_l) \ge (c-\epsilon)\frac{y_l}{\log^{2} y_l} & \text{   where   } & 
Q(X):=\{\P \in Q : \mathfrak{N}_{\K/\Q}(\P) \le X\}.
\end{eqnarray*}

Since every $\P \in Q$ has degree $1$, $\eta$ generates 
$(\mathcal{O}_{\K_{s}}/\mathfrak{r})^{\times}$ where $\mathfrak{r}
= \mathfrak{p}\cap \K_{s}$. Note $\mathfrak{N}_{\K/\Q}(\P)~=~p$ and $p\equiv u_1\bmod f$, $u_1 \bmod f$ restricts to a 
 generator of $C_s$ (by \lemref{generator}) and therefore to a generator
 of $Gal(\H(\K_s)/\K_s)$ and is trivial on $\K_s$.
 
 Therefore for all the primes in $\K_s$ below the ones in $Q$, the Artin symbol
 with respect to $Gal(\H(\K_s)/\K_s)$ corresponds to
 a generator of $C_s$. By 
 the isomorphism between said Galois group
 and the class group of $\K_s$ and by  \thmref{HarpCritvar}, we see that
this class must be a Euclidean ideal class.  This
completes the proof of \thmref{main-cubic}.

\section{Examples}
In this section, we provide some examples
as an application of \thmref{main-cubic}.
\begin{lem}\label{exlem}
Let $p$ and $q$ be two distinct primes which are $1 \bmod 3$.
The cyclotomic field $\Q(\zeta_{pq})$ contains four
subfields of degree $3$. If the class number of a
subfield not contained in both $\Q(\zeta_{p})$ and $\Q(\zeta_q)$ 
is $3$, then the Hilbert
class field of the subfield is abelian over $\Q$.
\end{lem}
\begin{proof}
The cyclotomic fields $\Q(\zeta_p)$
and $\Q(\zeta_q)$ have each a unique
subfield of degree $3$. Consider the compositum
of these fields. This is a subfield of $\Q(\zeta_{pq})$
of degree $9$ with Galois group $\Z/3\Z \times \Z/3\Z$.
This gives rise to four subfields of degree $3$, two of 
which are contained in $\Q(\zeta_p)$ and $\Q(\zeta_q)$.
Conversely any element of order $3$ in the Galois
group of $\Q(\zeta_{pq})/\Q$ must be in the unique subgroup isomorphic to $\Z/3\Z \times \Z/3\Z$
in the Galois group of of $\Q(\zeta_{pq})/\Q$
and therefore these are the only degree $3$ subfields.
For the two abelian cubic fields 
which are not contained in $\Q(\zeta_p)$ and $\Q(\zeta_q)$,
 the conductor is 
$pq$. Therefore both $p$ and $q$ ramify in these fields.
Now consider the following diagram.
$$
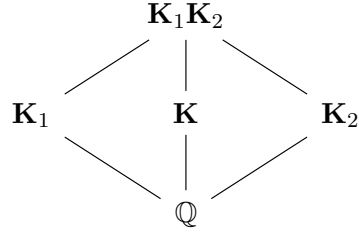
\begin{tikzcd}
&
\K_1\K_2 \arrow[dash]{dl}
\arrow[dash]{d}
\arrow[dash]{dr}\\
\K_1 \arrow[dash]{dr}{}
& 
\K
\arrow[dash]{d}
&
\K_2 \arrow[dash]{dl}{}\\
&
\Q
\end{tikzcd}
$$
Suppose that $\K_1 \subset \Q(\zeta_{p})$,
$\K_2 \subset \Q(\zeta_{q})$ and $\K$ is one of
the other two subfields. Then $p$
ramifies in $\K_1$ and not $\K_2$ and therefore
has ramification index $3$ with respect to $\K_1\K_2$.
Similarly $q$
ramifies in $\K_2$ and not $\K_1$ and therefore
has ramification index $3$ with respect to $\K_1\K_2$.
Since both $p$ and $q$ have ramification index $3$
with respect to $\K$, the extension $\K_1\K_2/\K$ is unramified 
at all places of $\K$. If the class number of $\K$
is $3$ then the Hilbert class field of $\K$ must be $\K_1\K_2$
and is therefore abelian over $\Q$.
\end{proof}

\begin{cor}\label{excorollary}
Consider four distinct primes $p_1,q_1,p_2$ and $q_2$,
all of which are $1 \bmod 3$. Suppose that
either $\K_1$ is a cubic
subfield of $\Q(\zeta_{p_1q_1})$ not contained
in $\Q(\zeta_{p_1})$ or $\Q(\zeta_{q_1})$
with class number
$3$ or it is a cubic subfield with class number $1$.
Similarly either $\K_2$ is a cubic subfield of 
$\Q(\zeta_{p_2q_2})$ not contained
in $\Q(\zeta_{p_2})$ or $\Q(\zeta_{q_2})$ 
with class number
$3$ or it is a cubic subfield with class number $1$.
Then at least one $\K_i$ must have a
Euclidean ideal class. 
\end{cor}

\begin{proof}
Since the conductors of $\K_1$
and $\K_2$ are co-prime to each
other, it follows that $\K_1\K_2$
over $\K_i$ must be ramified
for $i \in \{1,2\}$. 
By \lemref{exlem}, the Hilbert
class field of each $\K_i$ is abelian over $\Q$.
Therefore from \thmref{main-cubic}
we have the corollary.
\end{proof}

\begin{rmk}
We would like to remark here that
 if the Hilbert class field
is abelian over $\Q$
and the class number is a power
of $3$ greater than or equal to $9$,
 the class group
will no longer be cyclic (see \cite{MI}, Corollary to Theorem 5, Page 64-65).
\end{rmk}
The following table contains examples
of fields satisfying the hypothesis of
\corref{excorollary} (generated using
SAGE).

\begin{table}[ht]
\caption{Examples of cubic subfields with class number $1$ in $\Q(\zeta_{pq})$.}
\begin{tabular}{c c c c c}
\hline\hline
Serial & p & q & Defining polynomial of the subfield & Class number of the ring of integers \\
\hline
1 & 7 &  13 & $ x^3 - x^2 - 2x + 1$ & 1\\ 
2 & 7 & 13& $ x^3 - x^2 - 4x - 1 $ &1 \\
3 & 7 & 19& $ x^3 - x^2 - 2x + 1 $ &1 \\
4 & 7 & 19 & $ x^3 - x^2 - 6x + 7 $ &1 \\

\hline
\end{tabular}
\end{table}
\begin{table}[ht]
\begin{tabular}{c c c c c}
\hline\hline
Serial & p & q & Defining polynomial of the subfield & Class number of the ring of integers \\
\hline
5 & 7 &  31 & $ x^3 - x^2 - 2x + 1 $ &1 \\
6 & 7 & 31 & $ x^3 - x^2 - 10x + 8$ & 1 \\
7 & 7 & 37 & $ x^3 - x^2 - 2x + 1$ &1 \\
8 & 7 & 37 & $ x^3 - x^2 - 12x - 11$ &1\\ 
9 & 7 & 43 & $ x^3 - x^2 - 2x + 1 $ &1 \\
10 & 7 &  43 & $ x^3 - x^2 - 14x - 8$ & 1 \\
11 & 7 & 61 & $ x^3 - x^2 - 2x + 1 $ &1 \\
12 & 7 & 61 & $ x^3 - x^2 - 20x + 9 $ &1 \\
13 & 7 & 67 & $ x^3 - x^2 - 2x + 1$ & 1 \\
14 & 7 & 67 & $ x^3 - x^2 - 22x - 5 $ & 1 \\
15 & 13 & 19 & $ x^3 - x^2 - 6x + 7 $ & 1 \\
\hline
\end{tabular}
\end{table}
\begin{table}[ht]
\caption{Examples of cubic subfields with class number $3$ in $\Q(\zeta_{pq}).$\\
\phantom{mmmmmm}(not contained in $\Q(\zeta_p)$ or $\Q(\zeta_q)$)}

\centering
\begin{tabular}{c c c c c}
\hline\hline
Serial & p & q & Defining polynomial of the subfield & Class number of the ring of integers \\
\hline
 1 &7 & 13 & $ x^3 - x^2 - 30x - 27 $ & 3 \\
 2 &7 & 13 & $ x^3 - x^2 - 30x + 64 $ & 3 \\
 3 & 7 & 19 & $ x^3 - x^2 - 44x - 69 $ & 3 \\
 4 & 7 & 19 &  $ x^3 - x^2 - 44x + 64 $ & 3 \\ 
 5  & 7 & 31 & $ x^3 - x^2 - 72x - 209 $ & 3 \\
 6 & 7 & 31 & $ x^3 - x^2 - 72x + 225 $ & 3 \\ 
 7 & 7 & 37 & $ x^3 - x^2 - 86x + 211 $ & 3 \\
 8 & 7 & 37 & $ x^3 - x^2 - 86x - 48 $ & 3 \\ 
   9 & 7 & 43 &  $x^3 - x^2 - 100x + 379$ & 3 \\
10 & 7 & 43 &  $x^3 - x^2 - 100x - 223$ & 3 \\
11 & 7 & 61 & $x^3 - x^2 - 142x + 680$ & 3 \\
12 & 7 & 61 & $x^3 - x^2 - 142x - 601$ & 3 \\
13 &7 & 67 & $x^3 - x^2 - 156x - 608$ & 3 \\
14 &7  & 67 & $x^3 - x^2 - 156x + 799$ & 3 \\
15  & 13 & 19 & $ x^3 - x^2 - 82x + 64 $  & 3 \\
 \hline
\end{tabular}
\end{table}
\smallskip
\section{Data availability statement}
Data sharing is not applicable to this article as no datasets were generated or analysed during the current study.

\end{document}